# Non-parametric structural shape optimization of piecewise developable surfaces using discrete differential geometry

Makoto Ohsaki[1], Kentaro Hayakawa[1] and Jingyao Zhang[1]

[1] *Graduate School of Engineering, Kyoto University, Kyoto, Japan*

**Abstract** We propose a two-level structural optimization method for obtaining an approximate optimal shape of piecewise developable surface without specifying internal boundaries between surface patches. The condition for developability of a polyhedral surface onto a plane is formulated using the area of discrete Gauss map formed by unit normal vectors at the faces adjacent to each vertex. The objective function of the lower-level optimization problem is the sum of square errors for developability at all interior vertices. The contribution of large error to the objective function is underestimated by filtering with hyperbolic tangent function so that the internal boundary between the surface patches can naturally emerge as a result of optimization. Vertices are located non-periodically to generate the internal boundaries in various unspecified directions. Simulated annealing is used for the upper-level optimization problem for maximizing stiffness evaluated by the compliance under the specified vertical loads. The design variables are the heights of the specified points. It is shown in the numerical examples that the compliance values of the surfaces with a square and a rectangular plan are successfully reduced by the proposed method while keeping the developability of each surface patch. Thus, a new class of structural shape optimization problem of shell surfaces is proposed by limiting the feasible surface to piecewise developable surfaces which have desirable geometrical characteristics in view of fabrication and construction.

## 1. Introduction

Free-form architectural roof surfaces are generally designed using the continuous parametric surfaces such as Bézier surface and non-uniform rational B-spline (NURBS) surface [1]. The geometrical properties of surfaces such as normal vectors, principal curvatures, Gaussian curvature, and mean curvature are computed by complicated differentiation of the location vector of the surface, which is indirectly defined through the arc-length parameter, with respect to the parameters of the basis function. However, discretization into a polyhedral surface with triangular or quadrilateral mesh is necessary for evaluating structural performance using finite element analysis. Furthermore, for the latticed shells, a continuous surface should be approximated to triangular or quadrilateral mesh to cover the surface with planar or curved panels for designing an architectural roof. Therefore, it is desirable to use discretized model and formulations throughout the process of design, analysis, optimization, and construction. For this purpose, the theory and methods of discrete differential geometry [2] are effectively used.

Another important point for designing shell surfaces is that the cost for construction is drastically reduced if the surface consists of piecewise developable surface (PDS) patches that can be manufactured by bending plates without out-of-plane deformation or in-plane shear. However, the allowable surface shape is limited if the developable surface patch is modeled by parametric surfaces [3]. On the other hand, no such restriction exists if the surface patch is discretized into triangular or quadrilateral mesh. Hence, the formulations of discrete differential geometry can be used. Because the condition on a surface for developability onto the plane is formulated as the vanishing Gaussian curvature, the PDS is regarded as a special case of the surface with piecewise constant Gaussian curvature [4]. Using the formulations of discrete differential geometry, the Gaussian curvature and conditions for developability of the polyhedral surface onto the plane can easily be evaluated without carrying out complicated differential operations. The discrete formulation of developability condition is also utilized for designing a polyhedral surface using rigid origami [5].

One of the drawbacks of using polyhedral surfaces for generating developable surfaces is that the surface shape depends on the mesh discretization. Furthermore, the locations of internal boundaries between surface patches should be generally specified *a priori*. The authors developed a meshless approach, where only the locations of grid points are assigned to avoid mesh dependency [6, 7], and also proposed a method for underestimating a large error of developability to generate PDSs without specifying the locations of internal boundary [8]. However, structural optimization is not carried out in Ref. [8], and the locations of internal boundaries should be specified in Ref. [7].

In this study, we combine our methods in Refs. [7] and [8], and propose a two-level structural optimization method for approximately generating the stiffest PDSs without specifying internal boundaries between surface patches.

## 2. Two-level optimization problem

In this section, we briefly summarize the shape generation method proposed in Refs. [6–8] for completeness of the paper, and combine the shape generation problem of PDS and the structural shape optimization problem to propose a two-level optimization problem.

Consider a surface with grid points as shown in Fig. 1(a). Note that the dashed lines are drawn for visualizing the surface. Auxiliary edges adjacent to a point are generated as shown in Fig. 1(b). For example, point A is connected to eight neighborhood points with the red edges, and the normal vector is computed as the average of the unit normal vectors of eight triangular regions between the auxiliary edges.

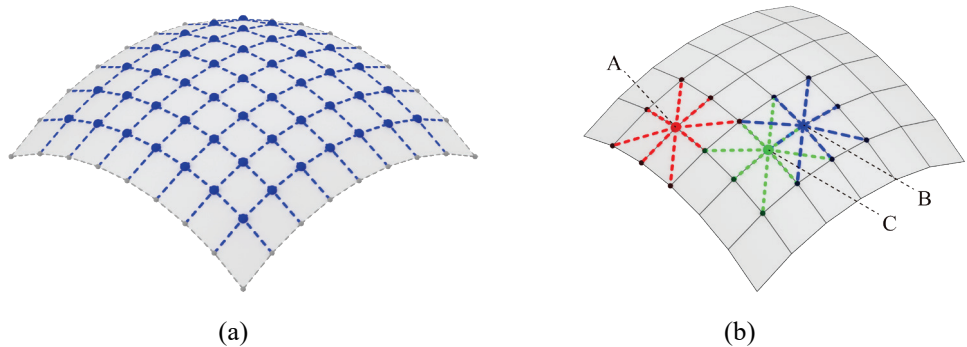

(a)                                (b)

Figure 1. Grid points and auxiliary edges on a curved surface [6].

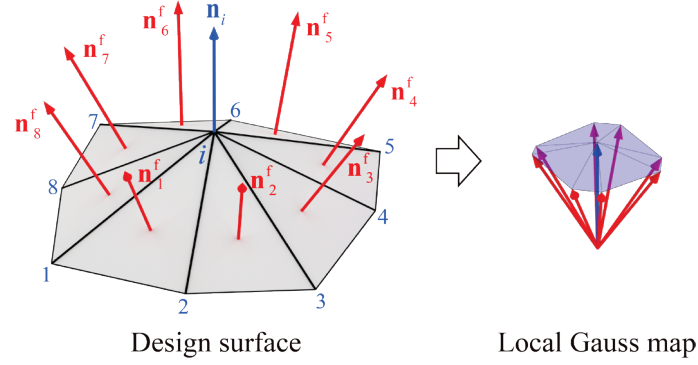

Figure 2. Local Gauss map at a point on a surface.

The developability condition at each grid point is formulated by the vanishing area of the Gauss map [6]. The unit normal vectors $\mathbf{n}_j^f$ ($j$ = 1, …, 8) of the triangular faces are defined as shown in Fig. 2, and the normal vector of grid point $i$ is computed as the average of the vectors on the adjacent eight faces. The unit normal vector at a grid point and those of adjacent faces are translated to have the same origin to generate the local Gauss map as shown in the right figure of Fig. 2. The area of local Gauss map vanishes if the surface is developable at the grid point.

A developable surface can be generated by minimizing the total area of Gauss map to zero. The vector of variables $\mathbf{x}$ consists of the coordinates of grid points. Let $a_{ij}$ denote the area of $j$th triangle of the Gauss map at grid point $i$. The sum of square of $a_{ij}$ of all adjacent triangles at node $i$ is denoted by $A_i$. In Refs. [6, 7], the sum of $A_i$ throughout all internal grid points in the set $I$ is minimized to generate a developable surface. In this case, the points for ignoring developability conditions, i.e., those to be excluded from $I$, should be assigned *a priori*, to generate PDSs. Therefore, in Ref. [8], the hyperbolic tangent function is used to underestimate a large error of developability condition so that the locations of internal boundaries need not be specified. Let $\mathbf{Z}$ denote the vertical coordinates of some selected points which are the vector of design variables of upper-level structural optimization problem. Then, the lower-level optimization problem for generating PDSs is formulated as

$$\text{Minimize} \quad F(\mathbf{x};\mathbf{Z}) = \sum_{i \in I} \tanh c(\sqrt{A_i(\mathbf{x};\mathbf{Z})+\varepsilon}) \quad (1)$$
$$\text{subject to} \quad \mathbf{x} \in \chi$$

Where $\chi$ is the feasible region of $\mathbf{x}$, and the shape of hyperbolic tangent function changes with the value of parameter $c$ as illustrated in Fig. 3. Note that $\mathbf{Z}$ is regarded as parameter vectors in Problem 1.

Let $W(\mathbf{Z})$ denote the compliance of the shell under specified loads. The upper-level structural optimization problem is formulated as

$$\text{Minimize} \quad W(\mathbf{Z}) = \tilde{W}(\mathbf{x}(\mathbf{Z}),\mathbf{Z}) \quad (2)$$
$$\text{subject to} \quad \mathbf{Z} \in \Xi$$

where $\Xi$ is the feasible region of $\mathbf{Z}$, $\tilde{W}$ is the compliance that is regarded as a function of both $\mathbf{Z}$ and $\mathbf{x}$, and $\mathbf{x}(\mathbf{Z})$ is obtained by solving the lower-level problem.

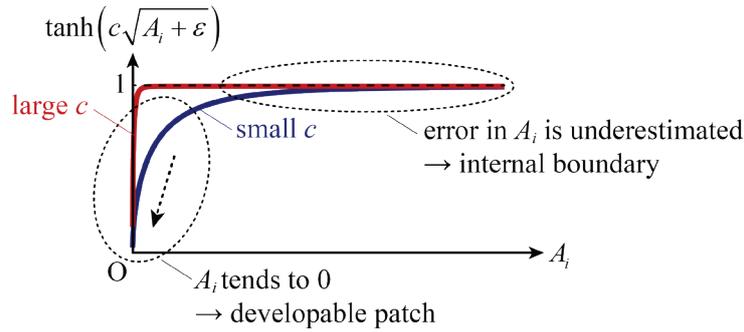

Figure 3. Conversion of developability error using hyperbolic tangent function.

## 3. Examples
3.1 Description of shell model and problem setting

We optimize two types of shells with square and rectangular plan, respectively, denoted by Cases 1 and 2. Problem 1 is solved using SLSQP [9] implemented in SciPy. The variables in Problem 1 are the $z$-coordinates of the points that are not fixed and are not included in the design variables of the upper-level problem. The lower and upper bounds of the $z$-coordinates at each grid point are $z_1 \pm 8$ (m) for the lower-level problem, and those for selected points as design variables are $z_2 \pm 1$ (m) for the upper-level problem, where $z_1$ and $z_2$ are the $z$-coordinate values at the initial shapes of Problems 1 and 2, respectively. Note that symmetry of the model is not considered, and the grid points are randomly moved in a horizontal direction so that the PDS obtained by solving Problem 1 does not strongly depend on the locations of grid points. The initial locations of the grid points are determined by randomly assigning the values of $u$ and $v$ parameters of the NURBS surface used for determining the entire shape of the surface. The parameters $u$ and $v$ are defined in the range from 0 to 1, and they are varied in the range of $\pm 0.015$ in Case 1 and in the range of $\pm 0.01$ in Case 2, respectively, from the equally spaced values.

Young's modulus and Poisson's ratio of the material are 20 GPa and 0.2, respectively, and the distributed load on the surface is 1.0 kN/m². Thickness of the shell is 0.1 m. Dual-annealing [10] implemented in SciPy is used for solving the upper-level structural optimization problem. The number of temperature reduction steps is 100, and local search is utilized. Neighborhood search is carried out 10 times at each temperature level; i.e., the total number of analysis is 1000. Finite element analysis is carried out using OpenSeesPy [11], and the quadrilateral shell element 'shellMITC4' is used. In the following examples, the right figures show contour of maximum bending moment computed as the larger value of the two principal bending moments at the center of element.

3.1 Case 1: Square plan

Optimal shapes are found from the initial shape in Fig. 4 that has a square plan. The size of shell on the plane is $10 \times 10$ m². The center point, four corner points, and the points at center of the four exterior boundaries are fixed in Problem 1. The $z$-coordinates of the five points which are fixed in Problem 1 but not at the corners as illustrated in Fig. 4(a) are the variables in the upper-level structural optimization. The boundary condition for structural analysis is given as shown in Fig. 4(a). The area of local Gauss map is plotted in Figs. 4(c), 5(c), and 6(c), where red collor corresponds to a large value indicating that developability condition is not satisfied. We can see from these figures that the points violating the developability condition are located around the center for the initial solution, and they spread to the points around the corners and the centers of edges for the optimal solution of Problem 2.

The value of compliance $W$ of the initial non-developable surface in Fig. 4 is 4.998 kNm, which increases to 5.555 kNm by obtaining a PDS in Fig. 5 by solving Problem 1. We can see from Figs. 4(b) and 5(b) that the maximum bending moment increases along the internal boundaries. The optimal shape of Problem 2 is shown in Fig. 6, where $W$ is reduced to 3.386 kNm by structural optimization. Note that irregularity exists in the shape near the fixed support and the points at which the $z$-coordinates are considered as design variables in Problem 2. The solutions for $c=10$ and 200 for defining the shape of hyperbolic tangent function are shown in Fig. 7. No clear internal boundary appears if $c$ has a small value 10, whereas a value of $c$ larger than 100 does not have strong effect on the optimal shape.

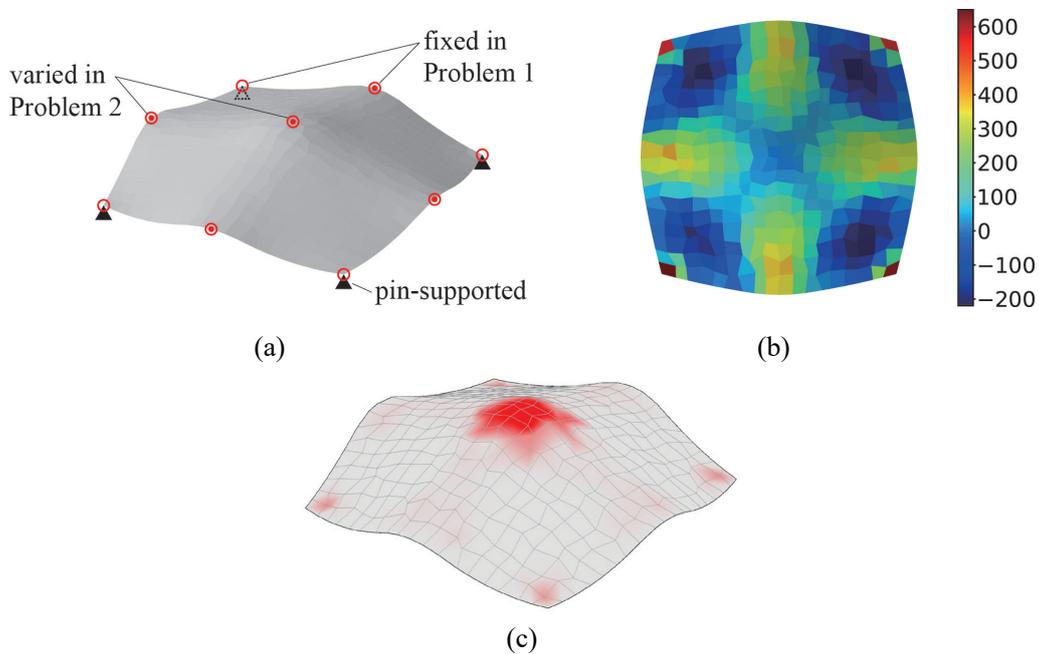

Figure 4. Initial shape of Problem 1 of initial solution of Problem 2 (Case 1);
(a) diagonal view, (b) maximum bending moment, (c) area of local Gauss map.

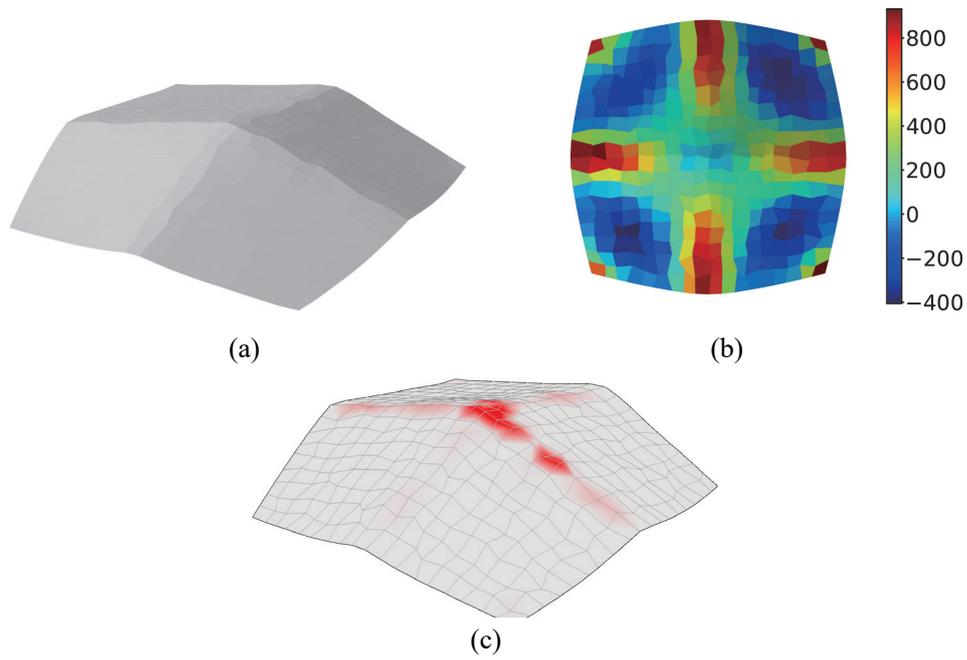

Figure 5. Optimal shape of Problem 1 of initial solution of Problem 2 (Case 1); $c = 100$;
(a) diagonal view, (b) maximum bending moment, (c) area of local Gauss map.

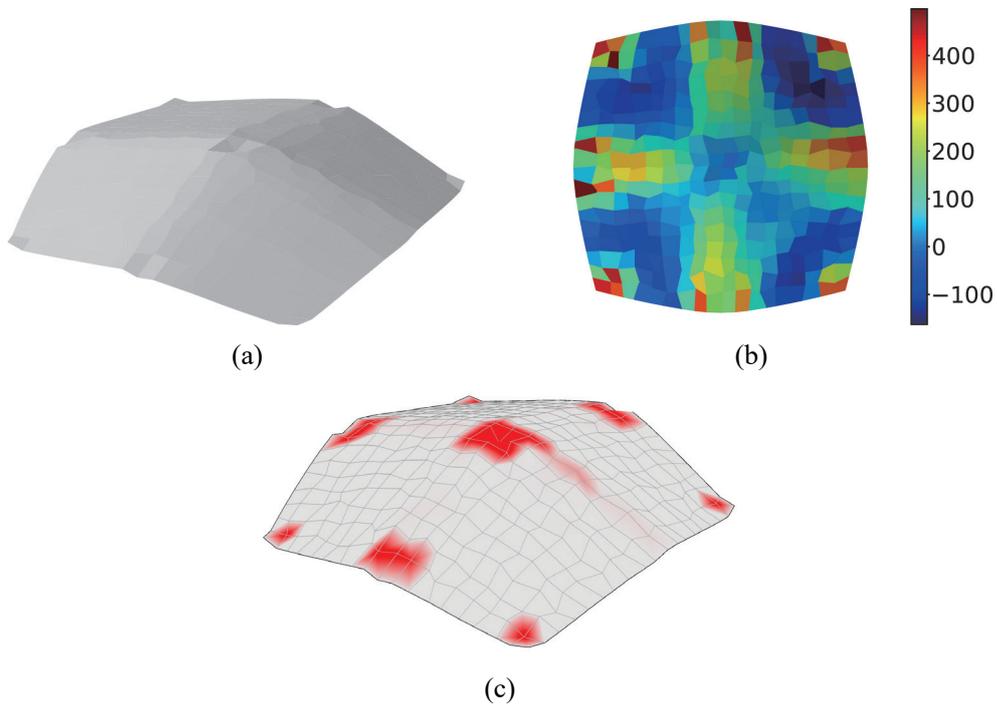

Figure 6. Optimal shape of Problem 2 (Case 1); $c = 100$;
(a) diagonal view, (b) maximum bending moment, (c) area of local Gauss map.

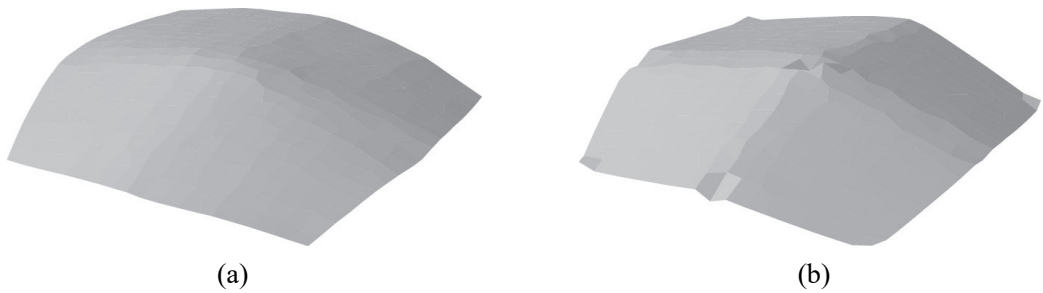

Figure 7. Optimal shapes of Problem 2 (Case 1); (a) $c = 10$, (b) $c = 200$.

3.2 Case 2: Rectangular plan

Optimal shapes are next found from the initial shape in Fig. 8 that has a rectangular plan. The size of shell on the plane is $10 \times 5$ m. Fifteen points shown in Fig. 8 are fixed in Problem 1. The $z$-coordinates of the eleven points which are fixed in Problem 1 but not at the corners are the variables in the upper-level structural optimization. The developability condition is not considered at the four red points. The boundary condition for structural analysis is given as shown in Fig. 8.

The value of compliance $W$ of the initial non-developable surface in Fig. 8 is 0.1480 kNm, which increases to 0.4338 kNm by obtaining a PDS in Fig. 9 by solving Problem 1 with $c = 10$. It is seen from Fig. 9(a) that the obtained PDS has a pair of conical surfaces connected along the curved internal boundary. We can see from Figs. 8(b) and 9(b) that the maximum bending moment has large values at the four corners, and those at the points along the internal boundary increase as a result of generating a PDS.

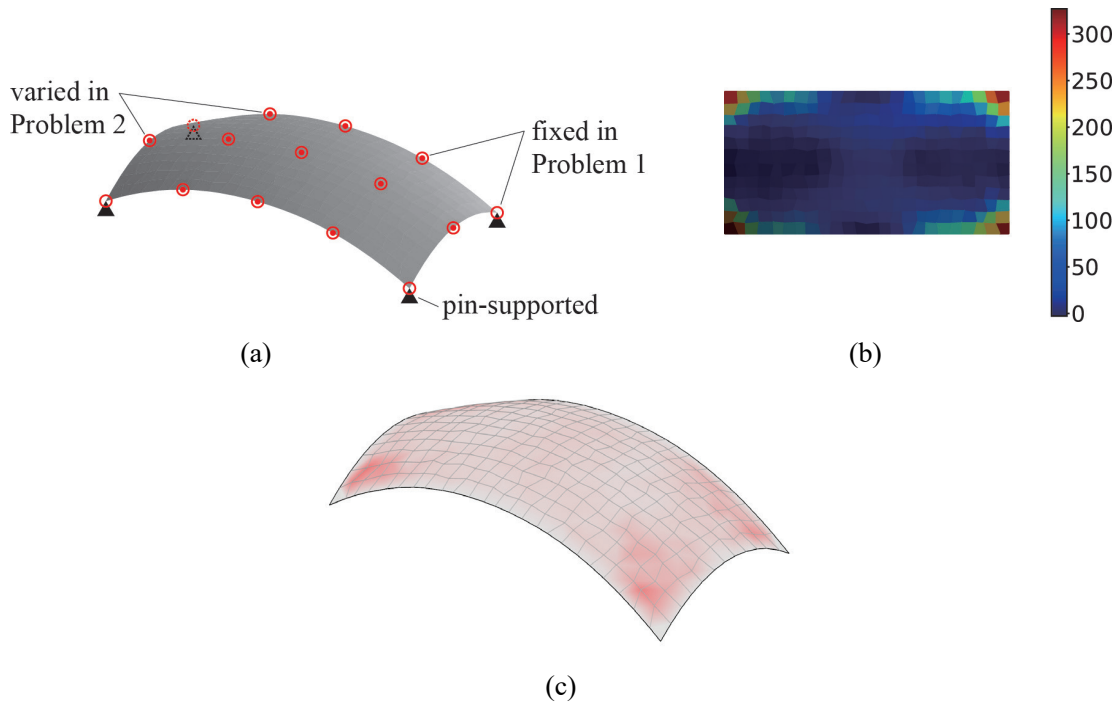

Figure 8. Initial shape of Problem 1 of initial solution of Problem 2 (Case 2);
(a) diagonal view, (b) maximum bending moment, (c) area of local Gauss map.

The optimal shape of Problem 2 is shown in Fig. 10, where $W$ is reduced to 0.06668 kNm by structural optimization. We can see from Fig. 10(a) that the optimal PDS has a cylinder, which is almost plane, around the center, and cylindrical surface patches are connected by ridges in the two side regions. The area of local Gauss map is plotted in Figs. 8(c), 9(c), and 10(c). We can confirm from Fig. 9(c) that the points violating the developability condition are located along the internal boundary of the two developable surfaces in Fig. 9. However, those points are not distinctly distributed for the optimal solution of Problem 2. The optimal shape for $c = 100$ is plotted in Fig. 11. Although the convergence near the variable points of problem 2 is not good, a surface consisting of eight planar patches is observed in Fig. 11.

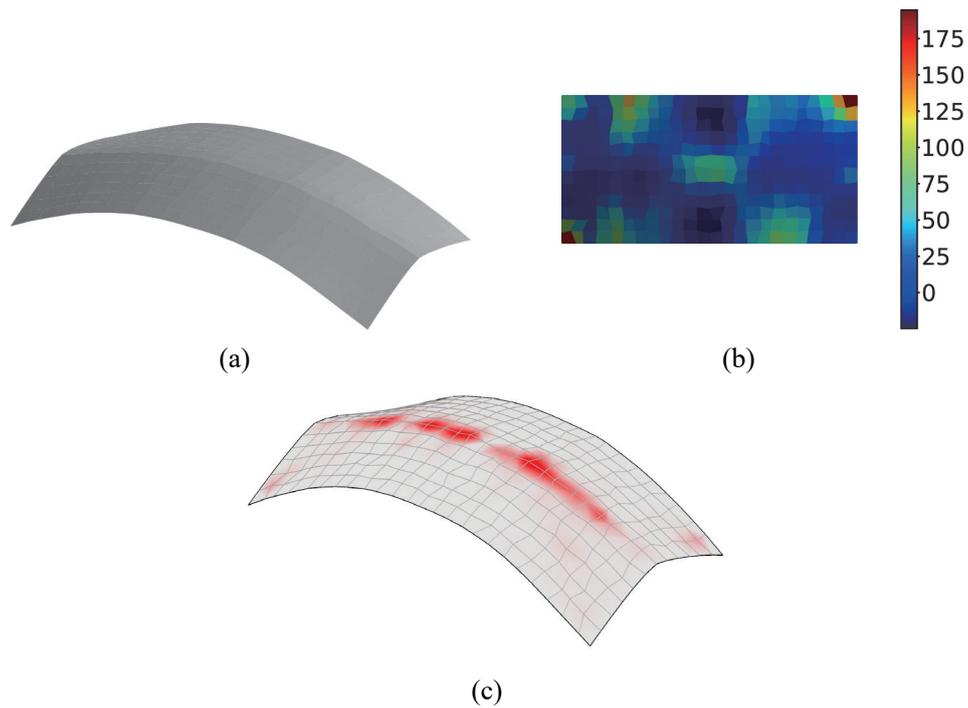

Figure 9. Optimal shape of Problem 1 of initial solution of Problem 2 (Case 2); $c = 10$;
(a) diagonal view, (b) maximum bending moment, (c) area of local Gauss map.

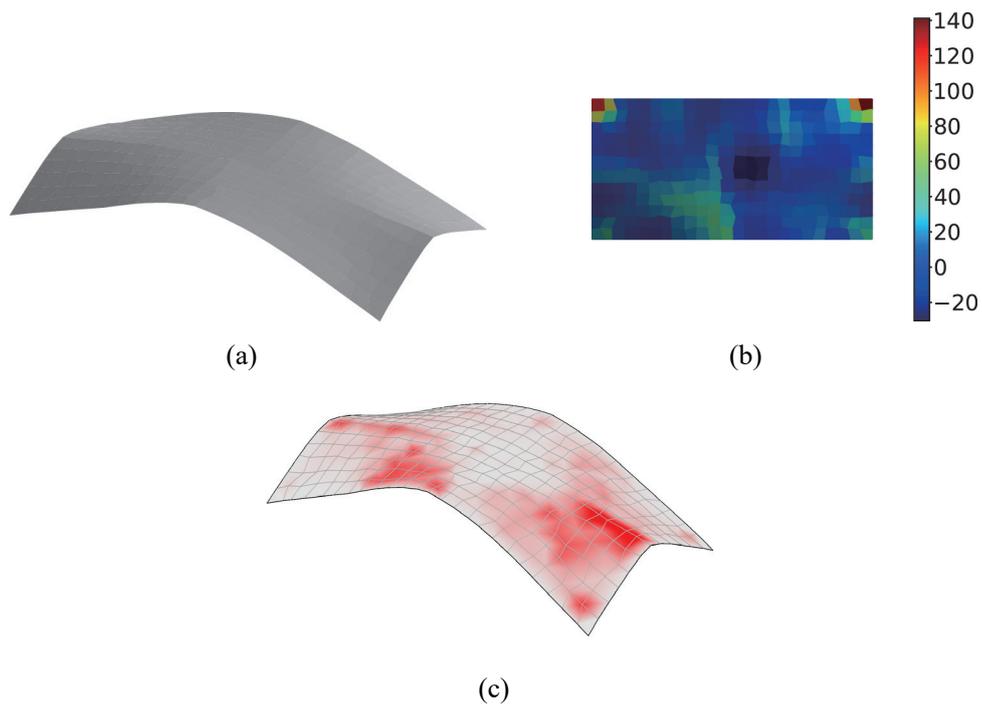

Figure 10. Optimal shape of Problem 2 (Case 2); $c = 10$;
(a) diagonal view, (b) maximum bending moment, (c) area of local Gauss map.

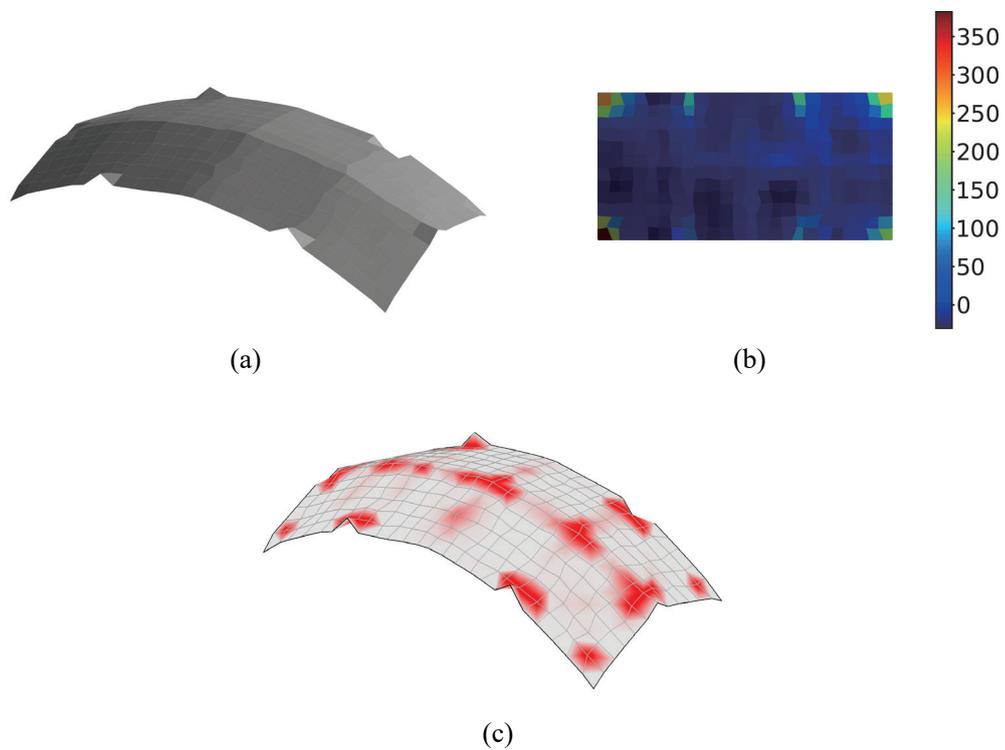

Figure 11. Optimal shape of Problem 2 (Case 2); $c = 100$;
(a) diagonal view, (b) maximum bending moment, (c) area of local Gauss map.

## 4. Conclusions

A two-level structural optimization method has been proposed for approximately generating the stiffest PDSs without specifying internal boundaries between surface patches.

It has been shown in the numerical examples that the values of objective function, which is the compliance under the specified load, of surfaces with square and rectangular plans are successfully reduced by the proposed method while keeping developability of each surface patch. Thus, a new class of shape optimization problem of shell surfaces has been proposed by limiting the feasible surface to PDSs which have desirable geometrical characteristics in view of fabrication and construction.

**Acknowledgment**

This work was partly supported by the Japan Science and Technology Agency, CREST Grant Number JPMJCR1911.